\documentclass[10pt]{amsart}

\usepackage{amsmath}
\usepackage{amsfonts}
\usepackage{amssymb}
\usepackage{amscd}
\usepackage{amsthm}
\usepackage{amsfonts}

\usepackage[all]{xy}
\usepackage{graphicx}
\usepackage{psfrag}

\usepackage{plain}
\usepackage[all]{xy}

\usepackage[active]{srcltx}
\usepackage{colordvi}

\newtheorem{prop}{\bf Proposition}[section]
\newtheorem{cor}[prop]{{\bf Corollary}}
\newtheorem{lem}[prop]{{\bf Lemma}}
\newtheorem{thm}[prop]{{\bf Theorem}}
\newtheorem{defn}[prop]{{\bf Definition}}

\numberwithin{equation}{section}

\newenvironment{rem}{{\bf Remark }}{ }

\newenvironment{pf}{{\bf Proof: }}{\qed\endtrivlist}

\newcommand{\Z}{\mathbb{Z} }

\newcommand{\C}{\mathbb{C} }
\newcommand{\R}{\mathbb{R} }

\newcommand{\ind}{\operatorname{ind}}
\newcommand{\Id}{\operatorname{Id}}

\newcommand{\Ch}{\operatorname{Ch}}

\newcommand{\supp}{\ensuremath{supp}}

\begin{document}

\title{A coarse Relative-partitioned index theorem}
\author{M. Karami, M.E. Zadeh, A. Sadegh}

\address{Mostafa Esfahani Zadeh,
Sharif University of Technology,
Tehran-IRAN
\newline 
Moin Karami, Sharif University of Technology, Tehran-Iran
\newline 
Ahmad R. Sadegh, Pennsylvania State University}
\email{esfahani@sharif.ir}

\begin{abstract}
It seems that the index theory for non-compact spaces has found its ultimate formulation 
in realm of coarse spaces and $K$-theory of related operator algebras. 
Relative and partitioned index theorems may be mentioned as two important and interesting examples 
of this program.      
In this paper we formulate a combination of these two theorems and establish a 
partitioned-relative index theorem.  
\end{abstract}

\maketitle

\section{Introduction}
For $i=1,2$ let $(M_i,g_i)$ be non-compact odd-dimensional spin complete Riemannian manifolds 
which are  partitioned by  hyper-surfaces $N_i\subset M_i$ into sub manifolds $M_i^+$ 
and $M_i^-$ with common boundary $\partial M_i^+=\partial M_i^-=M_i^+\cap M_i^-=N_i$. 
Let $E_i$ be a Clifford bundles over $M_i$ and put $H_i=L^2(M_i,E_i)$. 
The associated  Dirac operator $\mathbf D_i$ is a formally self adjoint 
operator on $H_i$. 
  The coarse index 
$\index{D_i}$ is an element in the $K$-theory of the coarse $C^*$-algebra $C^*(M_i)$ \cite{Roe-index}.
The restriction of operator $\mathbf D_i$ to $N_i$ is the grading reversing 
Dirac operator $D_i$ that acts on graded sections of $E_{|N_i}$.

Moreover assume there are closed subsets $W_i\subset M_i$ that intersect $N_i$ 
coarsely and transversally such that the intersection $Z_i=W_i\cap N_i$ is compact. 
We assume also that there is an isometry $\psi:M_1\verb+\+W_1\to M_2\verb+\+W_2$ 
that is covered by an isometry of bundles that we denote by $\Psi$, 
so that $\mathbf D_2=\Psi\circ\mathbf D_1\circ\Psi^{-1}$ on the sections supported in $M_2\verb+\+W_2$. Assume there is another 
Riemannian manifold $M$ and subsets $N$ and $W$ satisfying above conditions, 
(e.g. the intersection $Z:=N\cap W$ is compact and etc.) and there are smooth coarse maps 
$f_i\colon M_i\to M$ such that $f_i^{-1}(W)\subset W_i$ and  $f_i^{-1}(N)= N_i$. We 
assume that $N$ is a regular sub-manifold for $f_i$. 
Using all these structures (except the partitioning hyper-surfaces), we will define a 
relative index $\ind(\mathbf D_1,\mathbf D_2)\in K_1(C^*(W\subset M;A))$.

Under above conditions,  $Z_i$ is a compact subset of $N_i$ and all geometric data 
on $N_i\verb+\+Z_i$ are identified via isomorphism $\psi_{|N_1\verb+\+Z_1}$. Also we have the maps 
$f_i$ from $N_i$ into $N$ with  $f_i^{-1}(Z)\subset Z_i$. Therefore one can define the 
relative index $\ind(D_1,D_2)$ as an element in $K_0(\mathbf K)\simeq\Z$ 
(this because $Z$ is compact), where $\mathbf K$ is the algebra of all compact operators on 
a separable infinite dimensional Hilbert space. The relative index $\ind(D_1,D_2)$ was 
introduced by J. Roe. We introduce the notion of a $W$-coarse function as follows:

\begin{defn}
A function $h:M\to \R$  is  $W$-coarse if its restriction to each bounded neighbourhood  
of $W$ is a coarse map. 
\end{defn}

We will show (see \ref{ycoarse}) that $h$ induces a natural morphisms 

\[h_*:K_1(W\subset M)\to K_1(\R)\simeq\Z\]
Moreover, if $h$ is smooth with regular value $0$ and $N=h^{-1}(0)$ then we have the following 
equality (see theorem \ref{mainthm}) which is a 
main result of this paper
\[ h_*(\ind(\mathbf D_1,\mathbf D_2))=\ind(D_1,D_2)\]

We use this theorem to re-prove the non-existence of a metric on $\tilde N\times \R$ with 
uniformly positive scalar curvature, provided that $\tilde N$ is a compact enlargeable manifold. 
This theorem was first proved (amongst other important facts) in 
\cite{GrLa3}\footnote{Actually Gromov-Lawson have proved the stronger result that there is 
no Riemannian metric with positive scalar curvature on $\tilde N\times\R$. 
As far as we know there is no index-theoretic proof for this result.}. 
In section \ref{sec.2} and \ref{sec.3} we establish necessary tools from  $K$-theory 
and index theory of coarse spaces in relative context and we formulate the statement 
of the main theorem \ref{mainthm}. In section \ref{sec.4} we prove this theorem by 
reducing the problem to cylindrical case. In section \ref{sec.5} we provide an 
application to the main theorem. 
As a new result, we prove that if $\tilde N$ is a compact enlargeable manifold then 
$\tilde N\times\R^n$ does not admit a Riemannian metric with uniformly positive scalar curvature coarsely  
equivalent to a product metric $dx^2+g_{\tilde N}$. This is specially meaningful when is compared to a 
theorem by Gromove and Laswon stating that there is always a metric with uniformly positive scalar curvature
on $\tilde N\times \R^n$, for $n\geq 3$ \cite[page 298]{GrLa3}.

\section{Some constructions in relative K-theory of coarse spaces}\label{sec.2}

Let $X$ be a complete proper metric space and $H$ be a Hilbert space which 
is an ample module 
over $C_0(X)$. The module action of a function $\phi$ on $H$ is denoted by 
$\rho(\phi)$ 
or just by $\phi$ if there is no risk of confusion. 
A bounded linear operator $T$ on $H$ is called controlled (or have 
finite propagation property) 
if there is $r>0$ such that for any 
$\phi$ and $\psi$ in $C_0(X)$ the relation $d(\supp(\phi), \supp(\psi))>r$ implies 
$\rho(\phi)T\rho(\psi)=0$. 
The operator $T$ is called pseudolocal if $\rho(\phi) T\rho(\psi)$ is 
compact provided that $\phi\psi=0$. 
The operator $T$ is locally compact if for $\phi$ as in above, 
the linear maps $\phi T$ and 
$T\phi$ are compact operators. 
Given a closed subset $Y\subset X$, the operator $T$ is supported near $Y$ if there is a 
constant $r$ such that 
$\rho(\phi)T=0=T\rho(\phi)$ if $d(\supp(\phi), Y)>r$.
 
Using the above definitions, the following $C^*$-algebras are defined: 
The space of all bounded pseudolocal operator on $H$ is denoted by 
$\mathcal{D}^*(X)$ while  
$D^*(X)$ consists of all bounded, controlled and pseudolocal operators on $H$. 
The space of all bounded and locally compact operators on $H$ is 
denoted by $\mathcal{C}^*(X)$ and   
$C^*(X)$ consists of all bounded, controlled and locally compact operators 
on $H$. 
It is easy to verify that $C^*(X)$ is an ideal of $D^*(X)$. The relative $C^*$-algebra 
$C^*(Y\subset X)$ is the ideal of $C^*(X)$ (and of $D^*(X$)) consisting 
of those operators which 
are supported in a bounded neighbourhood of $Y\subset X$ (i.e., supported near $Y$). 
Similarly $D^*(Y\subset X)$ is the 
ideal of $D^*(X)$  consisting of those operators 
which are supported near to $Y\subset X$. It is easy to see 
that $C^*(Y\subset X)$ is an ideal 
of $D^*(Y\subset X)$. Similarly we can define the algebra 
$\mathcal D^*(Y\subset X)$ and its ideal 
$\mathcal C^*(Y\subset X)$

In these definitions we have not mentioned the Hilbert space $H$ because the $K$-theory of 
these algebras are canonically independent of $H$. 
When we need to emphasize the Hilbert space we include it in the notation, 
e.g., $C^*(X, H)$.

The $K$-homology $K_j(X)$ and relative $K$-homology $K_j(Y\subset X)$
are defined by following relations, cf. \cite{RoeSiegel-Sheaf}
\begin{gather*}
 K_j(X)=K_{j+1}(\mathcal D^*(X)/\mathcal C^*(X))~~;~~
 K_j(Y\subset X)=K_{j+1}(\mathcal D^*(Y\subset X)/\mathcal C^*(Y\subset X))
\end{gather*}

It turns out that the following equalities hold \cite[Page 6]{siegel2012mayer}

\begin{gather*}
 \mathcal D^*(X)/\mathcal C^*(X)= D^*(X)/ C^*(X)~\text{ and }~
 \mathcal D^*(Y\subset X)/\mathcal C^*(Y\subset X)= D^*(Y\subset X)/ C^*(Y\subset X)
\end{gather*}
Therefore the K-theory long exact sequences  
associated to the following short exact sequences
\begin{gather}
0\to C^*(X)\to D^*(X)\to D^*(X)/C^*(X)\to 0\label{alef0}\\
0\to C^*(Y\subset X)\to D^*(Y\subset X)\to D^*(Y\subset X)/C^*(Y\subset X)
\to 0\label{alef1}
\end{gather} 
and their naturality  provide the following commutative 
diagram with exact rows 

\begin{equation}\label{alef2}
\begin{CD}
\dots  @>>> K_{j+1}(D^*(Y\subset X)) @>>> K_{j}(Y\subset X) 
@>A>> K_j(C^*(Y\subset X)) @>>>\dots\\
@.                 @VVV                        @VVV               
@VVV                   @.\\
\dots  @>>> K_{j+1}(D^*(X)) @>>> K_{j}(X) @>A>> K_j(C^*(X)) @>>>\dots
    \end{CD}
\end{equation}
Here the vertical morphisms are induced by inclusion while the 
morphisms $A$'s are assembly maps(see \cite[p.11]{RoeSiegel-Sheaf}).

The closed subspace $Y$ is a complete metric space, so it has 
its own $K$-homology 
and coarse $C^*$-algebras. A suitable extension by zero defines a natural 
maps from these absolute objects to 
the corresponding relative objects associated to the inclusion 
$Y\subset X$. These inclusions 
define the following natural isomorphisms \cite[p.10]{siegel2012mayer}:
\begin{gather}
K_*(C^*(Y))\simeq K_*(C^*(Y\subset X))~;~
K_*(D^*(Y))\simeq K_*(D^*(Y\subset X))~;\label{relabs1}\\
K_*(Y)\simeq K_*(Y\subset X)\label{relabs2}
\end{gather}

Now we are going to state and prove the relative version of 
the Mayer-Vietoris exact sequence. For this purpose let $N\subset X$ be a 
closed subset that partitions $X$ into closed subset $X^+$ 
and $X^-$ such that $X^+\cap X^-=N$. We assume this partitioning be excisive, i.e. 
for any $r>0$, there is $r'>0$ such that 
\[N_r(X^+)\cap N_r(X^-)\subset N_{r'}(X^+\cap X^-)\] 
With this assumption  
the following Mayer-Vietoris exact sequence is stated and proved in \cite{YuMayer}
\begin{equation}\label{meyerabs}
K_{j+1}(C^*(X^+))\oplus K_{j+1}(C^*(X^-))\to K_{j+1}(C^*(X))
\xrightarrow{\partial_{mv}} K_{j}(C^*(N))\to K_{j}(C^*(X^+))\oplus K_{j}(C^*(X^-))
\end{equation}

Similar Mayer-Vietoris exact sequences for algebra $D^*$ 
and for the $K$-homology  
are established in \cite{siegel2012mayer}  and the naturality and commutativity of 
these exact sequences with 
respect to the lower row exact sequence in \eqref{alef2} is proved. As a result, the 
assembly map $A$ and the Mayer-Vietoris morphism commute with each other. 

We need to establish the Mayer-Vietoris exact sequence in the relative context. For 
this purpose, observe that 
the intersection $Z:=Y\cap N$ is a closed 
partitioning subset of  $Y$ making partition $Y=Y^+\cup Y^-$. 
We assume that the partition $Y=Y^+\cup Y^-$ is also excisive. 
 \begin{prop}\label{pro} 
 Under above assumptions, we have following commutative diagram, 
 where the rows are part of long exact sequences and vertical maps are assembly maps
 
 \begin{equation*}
    \begin{CD}
K_{j}(Y\subset X) @>>{\partial_{mv}}> K_{j-1}(Z\subset N) @>>> 
K_{j-1}(Y^+\subset X^+)\oplus K_{j-1}(Y^-\subset X^-) \\
      @V {A} VV              @V{A}VV @V{A}VV  \\
 K_{j}(C^*(Y\subset X))@>>{\partial_{mv}}> K_{j-1}(C^*(Z\subset N))  @>>> 
K_{j-1}(C^*(Y^+\subset X^+))\oplus K_{j-1}(C^*(Y^-\subset X^-)) 
\end{CD}
 \end{equation*}
\end{prop}
\begin{pf}
Actually this follows from a more general Mayer-Vietoris exact sequence 
for $C^*$-algebras.  Consider a triple of $C^*$-algebras $(A, I_1,I_2)$, where $I_1$ and 
$I_2$ are ideals in $A$ such that $A=I_1+I_2$. Then the following Mayer-Vietoris 
exact sequence holds and is natural, c.f.  \cite{YuMayer}
\begin{equation}\label{mayergen}
\to K_{j+1}(I_1)\oplus K_{j+1}(I_2)\to K_{j+1}(A)
\xrightarrow{\partial_{mv}} K_{j}(I_1\cap I_2)\to K_{j}(I_1)\oplus K_{j}(I_2)\to
\end{equation}
It turns out that the following equalities (and similar equalities for 
$D^*$ and $D^*/C^*$) hold

\begin{gather*}\label{relgen}
C^*(Y\subset X)=C^*(Y^+\subset X)+
C^*(Y^-\subset X)\\
C^*(Y^+\subset X)\cap C^*(Y^-\subset X))=C^*((Z\subset N)\subset (Y\subset X))
\end{gather*}

In view of these equalities we can apply \eqref{mayergen} and this leads to the desired 
commutative diagram. 
\end{pf}

Following \cite{SchickZadeh-Large}, a metric space $Y$ is flasque if it has an isometry 
$\phi$ which is homotopic to 
identity such that $\phi^k(Y)$ leaves any compact subset of $Y$ for sufficiently large $k$. An example of 
such spaces is $\R\times Z$ with metric $dt^2+g_Z$, where $Z$ is  manifold with Riemannian 
metric $g_Z$. 
If $Y$ is flasque then $K_*(C^*(Y))$, $K_*(D^*(Y))$ and $K_*(Y)$ vanish. 
Therefore by 
\eqref{relabs1} and \eqref{relabs2} we have also 
\begin{gather}\label{flasque}
K_*(C^*(Y\subset X))=K_*(D^*(Y\subset X))=K_*(Y\subset X)=0
\end{gather} 
Using these vanishing results and the proposition \ref{pro} we get the following corollary 

\begin{cor}\label{meyerassem}
 With the notation of proposition \ref{pro} if $Y^+$ and $Y^-$ are flasque, and the partition 
 $Y=Y^+\cup Y^-$ is execive,  
 then we get the following  commutative diagram where the horizontal arrows are isomorphisms
 \begin{equation*}
    \begin{CD}
K_j(Y\subset X) @>{\delta_{mv}}>> K_{j-1}(Z\subset N)\\
     @V\mathbf A VV @V\mathbf A VV           \\
K_{j}(C^*(Y\subset X) @>{\delta_{mv}}>> K_{j-1}(C^*(Z\subset N)
 \end{CD}
  \end{equation*}
\end{cor}
Now we study the functoriality of these (relative) $K$-groups. 
Let $X$ and $X'$ be smooth Riemannian manifolds with hermitian bundles $E$ and $E'$ and let 
$H=L^2(X,E)$ and $H'=L^2(X', E')$.
A continuous coarse function $f:X\to X'$ induces 
canonical morphism $f_*:K_*(C^*(X,H))\to K_*(C^*(X',H'))$,   
 $f_*:K_*(D^*(X,H))\to K_*(D^*(X',H'))$ and  $f_*:K_*(X)\to K_*(X')$. We need to extend these 
morphisms to relative cases, so we describe very briefly their constructions. 
An isometry $V:H\to H'$ covers $f$ topologically if for every compactly supported 
$\phi\in C_0(X')$, the operator 
$V^*\phi\,V-\phi\circ f$ is a compact operator on $H$. 
It covers $f$ coarsely and has propagation speed less 
than $\epsilon>0$ if for any compactly supported section $\xi\in H$ one has 
$\phi\,V(\xi)=0$ provided that $d'(\supp (\phi), f(\supp(\xi)))\geq\epsilon$. 
Given a continuous and coarse function $f$, one can construct the isometry $V$ that covers 
$f$ both topologically and coarsely (with arbitrary small propagation speed) as follows. 
Let $\{U_i\}_i$ be a Borel covering for $X'$ such that the diameter of each $U_i$ 
is less than $\epsilon/2$ and the intersections $U_i\cap U_j$ have measure zero. 
There is an isometry $T_i:L^2(f^{-1}(U_i), E)\to L^2(U_i, E')$ that topologically covers 
$f_{|f^{-1}(U_i)}$ \cite[lemma 5.2.4]{HigRoe-Analytic}. Then $V=\oplus_i T_i$ is an isometry 
that covers $f$ both topologically and coarsely with propagation speed less than $\epsilon$.  
It turns out that the morphism 
\begin{equation}\label{adj}
Ad_V:B(H)\to B(H')~~;~~Ad_V(T)=VTV^*
\end{equation}
restricts to morphisms from $C^*(X)$ to $C^*(X')$, and from $D^*(X)$ to $D^*(X')$. Therefore 
it gives rise to a morphism from $D^*(X)/C^*(X)$ to $D^*(X)/C^*(X')$. 
The induced morphisms on $K$-theory  are actually independents of the covering 
isometry $V$ that is used in their definitions. 

Now let $Y\subset X$ and $Y'\subset X'$. We call a map $h:(Y,X)\to (Y',X')$ 
(this means $h(Y)\subset Y'$) 
 $Y$-coarse if for each $r>0$ the restriction of $h$ to $N_r(Y)$ is a coarse map. 
\begin{lem}\label{ycoarse}
A map $h:(Y,X)\to (Y',X')$ which is continuous and $Y$-coarse induces canonical morphisms 
\begin{gather*}
h_*:K_*(C^*(Y\subset X))\to K_*(C^*(Y'\subset X'))\\
h_*:K_*(D^*(Y\subset X))\to K_*(D^*(Y'\subset X'))\\
h_*:K_*(Y\subset X)\to K_*(Y'\subset X')
\end{gather*}
\end{lem}
\begin{pf}
Let $V$ be the isometry from $H$ to $H'$ constructed in the above discussion that covers $f$ 
 topologically and coarsely and has propagation speed less than $\epsilon$. 
Since $f(Y)\subset Y'$ and $f$ is coarse in any bounded neighbourhood of $Y$ it follows that 
for each $s>0$ there is $s'>0$ such that $f(N_s(y))\subset N_{s'}(Y')$ 
(this is because $f$ is uniformly expansive in $N_s(Y)$). Let $T\in B(H)$ be supported in a
 finite distance of $Y$. Because $V$ has finite propagation, $Ad_V(T)$ has support in a finite 
 distance of $Y$. By continuity we conclude that 
 $Ad_V$ is a morphisms from $C^*(Y\subset X)$ to $C^*(Y'\subset X')$, 
 and from $D^*(Y\subset X)$ to $D^*(Y'\subset X')$. Therefore, it
gives rise to a morphism from $D^*(Y\subset X)/C^*(Y\subset X)$ to 
$D^*(Y'\subset X)/C^*(Y'\subset X')$. The induced morphisms in $K$-theory level 
are the desired morphisms.  
\end{pf}

\begin{rem}
Although in this section we have considered the general proper metric spaces, 
for our purposes in this paper we need just complete Riemannian manifolds or complete sub 
manifolds of such manifolds. 
In this case a rich index theory with values in the $K$-theory of the 
coarse $C^*$-algebras has been developed that we are going to introduce, 
while we keep our discussion as brief as possible.   
\end{rem}

\section{Some constructions in coarse index theory}\label{sec.3}

In this section we assume the geometric context and notation of the introduction. Namely for $i=1,2$ we 
have $(M_i, g_i)$, the subsets $W_i$ and partitioning hyper-surfaces $N_i$ and 
Dirac type operators $\mathbf D_i$ that act on smooth sections of Clifford bundles $E_i$. 
An essential part of our assumption is the isometry 
$\psi:M_1\verb+\+ W_1\to M_2\verb+\+ W_2$ that lifts to an isometry $\Psi$ of 
bundles and conjugates the Dirac type operators: 
\[\Psi^{-1}\mathbf D_2\Psi=\mathbf D_1 \text{ on } C^\infty(M_1\verb+\+ W_1, E_1)\]

We review very briefly the construction of relative indices following \cite{Roe-Positive}. 
Let $\chi$ be a normalization function, i.e., a continuous function on $R$ that goes toward to $\pm 1$ 
when $x\in \R$ goes toward $\pm\infty$. It turns out that $\chi(\mathbf D_i)\in D^*(M_i)$ 
while $\chi^2(\mathbf D_i)-I\in C^*(M_i)$. Therefore, if $\alpha=(1+\chi)/2$, 
then $\alpha(\mathbf D_i)$ is an idempotent 
in the quotient algebra $D^*(M_i)/C^*(M_i)$ and determines an element $[\alpha(\mathbf D_i)]$ in 
$K_0(D^*(M_i)/C^*(M_i))$. The short exact sequence  
\[0\to C^*(M_i)\to D^*(M_i)\to D^*(M_i)/C^*(M_i)\to 0\]
gives rise to a long exact sequence in $K$-theory: 
\[\dots\to K_{j+1}(D^*(M_i))\to K_{j+1}(D^*(M_i)/C^*(M_i))=K_j(M_i)\xrightarrow{\delta} 
K_{j}(C^*(M_i))\to\dots.\]
The coarse index of $\mathbf D_i$ is given by 
\[\ind(\mathbf D_i)=\delta([\alpha(\mathbf D_i)])\in K_1(C^*(M_i))\]

Since $N_i$ is even-dimensional, the restriction of $E_i$ to $N_i$ 
(that we denote by the same symbol $E_i$) is graded 
$E_i=E_i^+\oplus E_i^-$. Let $\theta:E_i^+\to E_i^-$ be a Borel bundle isomorphism and 
let $\chi$ be an odd normalization function and put $\beta=\theta^{-1}\chi$. 
Then $\beta(D_i^+)=\theta^{-1}\chi(D_i^+)$ belong to 
$D^*(N_i)$, while  $\beta^2(D_i^+)-I$; as an operator on $L^2(N_i, E_i^+)$; 
belongs to $C^*(N_i)$. Therefore 
$\beta(D_i^+)$ determines an element in $K_1(D^*(N_i)/C^*(N_i))$ and 
the index map $\delta$ defines the coarse index 
\[\ind(D_i)=\delta([\beta(D_i)])\in K_0(C^*(N_i))\]
Because $M_i$ (and $N_i$) are not compact, the indices $\ind(\mathbf D_i)$ 
(and $\ind(D_i)$) might not  be  
interesting. However, using $\psi$, it is possible to extract a 
part of these indices that reflects some geometric and topological 
contents of the difference of the operators and underlying spaces. This construction is the relative index and is 
introduced by J. Roe.  We follow his arguments as stated in \cite{Roe-Positive} 
and \cite{Roe-relative} to construct a partitioned version of  relative indices.

We give the definition of $\ind(\mathbf D_1,\mathbf D_2)$, the definition of 
$\ind(D_1, D_2)$ is similar and we will briefly discuss it. 
We recall from introduction the functions 
$f_i:M_i\to M$ with $W_i=f_i^{-1}(W)$. 
Let $f_i$ be continuous and coarse 
and let $V_i$ be an isometry between 
$H_i:=L^2(M_i, E_i)$ and $H:=L^2(M)$ that covers topologically and coarsely the map $f_i$. 
As we have mentioned earlier, for a given $\epsilon>0$ we may assume that the 
propagation speed of $V_i$ is at most $\epsilon/2$. Moreover we assume that $V_i$ maps 
$L^2(M_i\verb+\+W_i)$ on $L^2(M\verb+\+W)$ and $V_2\Psi=V_1$. 

Given $T_i\in B(H_i)$ then $Ad_{V_i}(T_i)=V_iT_iV_i^*\in B(H)$ and $Ad_{V_i}$ 
provides $C^*$-isomorphism from 
$C^*(M_i,H_i)$ (resp. $C^*(Y_i\subset M_i,H_i)$) to  
$C^*(M,H)$ (resp. $C^*(Y\subset M,H)$). Using these isomorphisms, we will 
consider in below $T_i$ as an element in $C^*(M)$. This discussion is also true for 
$D^*$ and for $K$-homology because $f_i$'s are continuous. 
With this discussion let $T_1$ and $T_2$ be respectively elements in $C^*(M_1)$ and 
$C^*(M_2)$ which are considered as elements in $C^*(M)$. We call them conjugate via 
$\psi$ up to $C^*(W\subset M)$, and denote it by $T_1\sim_\psi T_2$, 
if $\Psi^{-1}T_2\Psi =T_1 $  outside a bounded neighborhood of $W$ or 
if $(T_1, T_2)$ is the norm limit of such pairs of operators.  
Now consider following $C^*$-algebras 
\begin{gather*}
 \mathcal A=\{(T_1,T_2)|T_1\sim_\psi T_2\}\subset C^*(M)\times C^*(M)\\
\mathcal B=\{(T_1,T_2)| T_1\sim_\psi T_2\}\subset D^*(M)\times D^*(M) 
\end{gather*}
It is clear that $\mathcal A$ is an ideal in $\mathcal B$; therefore, one has 
the long exact sequence 
\begin{gather}\label{relcnc}
\dots\to K_0(\mathcal B)\to K_0(\mathcal B/\mathcal A)\xrightarrow{\delta} K_1(\mathcal A)\to\dots
\end{gather}

By above discussion, $(\alpha(\mathbf D_1), \alpha(\mathbf D_2))$ is an element 
in $\mathcal B$ and determines the class $[(\alpha(\mathbf D_1), \alpha(\mathbf D_2))]$ in 
$K_0(\mathcal B/\mathcal A)$, therefore 
$\delta([(\alpha(\mathbf D_1), \alpha(\mathbf D_2))])$ belongs to $K_1(\mathcal A)$. 
The following unitary map equals $\Psi$ on second summand
\begin{gather}\label{unit1}
U:=V_2^{-1}V_1:L^2(W_1,E_1)\oplus L^2(M_1-W_1, E_1)
\to L^2(W_2,E_2)\oplus L^2(M_2-W_2, E_2)
\end{gather}

Now consider the following short exact sequence 
\begin{gather}\label{short1} 
0\to C^*(W\subset M, H_1)\to \mathcal A\xrightarrow{\pi} C^*(M, H_2)\to 0
\end{gather}
where the first morphism is $a\to (a,0)$ and the second one is $(a,b)\to b$. 
This sequence splits and the splitting morphism is given by 
$(U^{-1}b\,U,b)\xleftarrow{q} b$, therefore 
\begin{gather}\label{dirs} 
K_1(\mathcal A)=K_1(C^*(W\subset M))\oplus K_1(C^*(M))
\end{gather}
We denote the projection on the first summand by $p$ (which is equal to $\Id- q\circ \pi_*$). 
The relative index is defined by 
\[\ind(\mathbf D_1,\mathbf D_2)=p\circ\delta([(\alpha(\mathbf D_1), \alpha(\mathbf D_2))])
\in K_1(C^*(W\subset M))~.\]
\begin{rem}
Similarly one defines the algebras 
\begin{gather*}
 A=\{(T_1,T_2)|T_1\sim_\psi T_2\}\subset C^*(N)\times C^*(N)\\
B=\{(T_1,T_2)|T_1\sim_\psi T_2\}\subset D^*(N)\times D^*(N) 
\end{gather*}
and the unitary $U$ 
\[U: L^2(N_1,E_1)=L^2(Z_1,E_1)\oplus L^2(N_1-Z_1, E_1)\to  
L^2(N_2,E_2)=L^2(Z_2, E_2)\oplus L^2(N_2-Z_2, E_2)\]
One has the direct sum relation 
$K_0(A)=K_0(C^*(Z\subset N))\oplus K_0(C^*(M))$. 
The relative index $\ind(D_1,D_2)$ is defined by a similar procedure 
but differs in parities: 
\[\ind(D_1,D_2)=p\circ\delta([(\beta(D_1), \beta(D_2))])\in K_0(C^*(Z\subset N))\simeq\Z~.\]
\end{rem}

We need also to define the relative class $[\mathbf D_1, \mathbf D_2]$ 
as an element in the $K$-homology group $K_1(W\subset M)$ and investigate its relation 
with relative index. For this purpose the first step is to note that short exact sequences 
similar to \eqref{short1} holds also for algebras $\mathcal B$ and $\mathcal B/\mathcal A$, i.e. 
the following sequences are exact and splits  

\begin{gather*}
0\to D^*(W\subset M, H_1)\to \mathcal B\xrightarrow{\pi} D^*(M, H_2)\to 0\\
0\to D^*(W\subset M, H_1)/C^*(W\subset M, H_1)\to \mathcal B/\mathcal A\xrightarrow{\pi} 
D^*(M, H_2)/C^*(M, H_2)\to 0
\end{gather*}
Each term in exact sequence \eqref{short1} is an ideal in corresponding term in the 
first exact sequence in above and the inclusion commute with other arrows and splitting 
morphisms. Therefore we can pass to long $K$-theory exact sequence to get 
the following commutative diagram with horizontal isomorphisms
\begin{equation}
\xymatrix{
K_*(\mathcal B/\mathcal A) \ar[d]^{A} \ar[r]
&K_*(W\subset M)\ar[d]^{A}& \oplus & K_*(M)\ar[d]^{A} \\
K_{*+1}(\mathcal A) \ar[r]& K_{*+1}(C^*(W\subset M)) &\oplus & K_{*+1}(C^*(M))}
\end{equation}
As we mentioned earlier $[\alpha(\mathbf D_1), \alpha(\mathbf D_2)]$ belongs to 
$K_0(\mathcal B/\mathcal A)$. We call the image of this class in $K_1(W\subset M)$ the relative $K$-homology class and denote it by 
$[\mathbf D_1, \mathbf D_2]\in K_1(W\subset M)$ (compare to \cite[p.15]{siegel2012mayer}). 
It is clear by this definition and the commutativity of above diagram that 
\begin{lem}\label{asemrel}
Let $A:K_0(W\subset M)\to K_1(C^*(W\subset M))$ be the assembly map; then 
\[A([\mathbf D_1, \mathbf D_2])=\ind(\mathbf D_1, \mathbf D_2)\]
\end{lem}

We recall that $N$ partitions $M$ into subsets $M^+$ and $M^-$.   
Let $h:M\to \R$ be the signed distance function from $N$ that is non-negative on 
$M^+$ and non-positive on $M^-$. Then $0$ is a regular value with $N=h^{-1}(0)$.  
\begin{lem}
$h$ is a continuous and $W$-coarse function.
\end{lem}
\begin{pf}
It is easy to see that $h$ is Lipschitz, so it is continuous and uniformly contractive. 
We need just to show that the restriction of $h$ to any bounded neighbourhood of $W$ is proper.   
Let $N_r(W)$ be a finite distance neighbourhood of $W$. For a given $r'>0$ we 
have $N_{r'}(N)=h^{-1}([0,r'])$. Because 
$N$ and $W$ intersects coarsely with intersection $Z$, there is $s>0$ such that 
$N_r(W)\cap N_{r'}(N)\subset N_s(Z)$. This last set is compact because $Z$ is compact. 
Therefore the restriction of $h$ to $N_r(W)\subset M$ is proper. 
\end{pf}
Consider the following isomorphisms where $P:N\to pt.$ is the constant map to 
a single point (here the point $0\in \R$) and $\partial_{mv}$ is the Mayer-Vietoris 
isomorphism corresponding to $\R=\R^{\geq0}\cup \R^{\leq0}$
\begin{gather}
K_1(C^*(\R))\xrightarrow{\partial_{mv}} K_0(C^*(pt.))=K_0(\mathbf K)\simeq \Z\label{conv1}\\
K_0(C^*(Z\subset N))\xrightarrow{P_*}K_0(C^*(pt.))=K_0(\mathbf K)\simeq \Z\label{conv2}
\end{gather}
Using above lemma and lemma \ref{ycoarse}, $h$ induces the following morphism 
\begin{equation}\label{naghol}
h_*:K_1(C^*(W\subset M))\to K_1(C^*(\R))
\end{equation}   
Now we can state the main theorem of this paper 
\begin{thm}\label{mainthm}
The following equality holds when each side is considered as an element in $\Z$ according to 
above isomorphisms
\begin{equation}\label{mainequ}
h_*(\ind(\mathbf D_1,\mathbf D_2))=\ind(D_1,D_2)
\end{equation}
\end{thm}
This theorem is actually true for $h$ a $W$-coarse map having $0$ as a regular value 
and $h^{-1}(0)=N$. 
 strategy to prove this theorem is reducing the general case to product 
case and then proving the product case. In reducing to product case we follow the 
approach of \cite{SchickZadeh-Large}, while to prove the product case we follow the approach 
of \cite{siegel2012mayer}.

\section{Proof of the main theorem}\label{sec.4}

The first step in giving a proof for the main theorem \ref{mainthm} is 
the following proposition.
\begin{prop}\label{redtopro}
The partitioned-relative index $h_*(\ind(\mathbf D_1,\mathbf D_2))\in K_1(\R)$ 
depends only on geometric 
data in  bounded neighbourhoods around $N_1\subset M_1$ and  $N_2\subset M_2$ .
\end{prop}
\begin{pf}
The proof of this proposition is based on a strong localization property for $K_*(\R^p)$ 
which is formulated and proved in \cite[Proposition 3.4]{SchickZadeh-Large}: let $T_1$ and $T_2$ be 
elements in $D^*(\R^p)$ which determine classes $[T_1]$ and $[T_2]$ in 
$K_*(\R^p)=K_{*+1}(D^*(\R^p)/C^*(\R^p))$. The localization property asserts 
that these classes are equal if $T_1=T_2$ in an open subset of $\R^p$. 
The naturality of the assembly maps and lemma \ref{asemrel} imply 
\[h_*(\ind(\mathbf D_1,\mathbf D_2))=h_*\circ A([\mathbf D_1,\mathbf D_2])
=A\circ h_*([\mathbf D_1,\mathbf D_2])\]
Therefore, due to the localization property, it is enough to show that the geometry of 
$\mathbf D_1$ and $\mathbf D_2$ around $N_1$ 
and $N_2$ determines the value of a representative of 
$h_*([\mathbf D_1,\mathbf D_2])\in K_1(\R)$ in a neighborhood of $0\in \R$. 
To do this we fix the compact neighborhood $K=N_{2a}(Z)$ of $Z$ in $M$. 
We recall that $0$ is a regular value of $h$ with $N=h^{-1}(0)$ and $N$ consists of 
regular values 
for $f_i$ for $i=1,2$ with $N_i=f_i^{-1}(N)$. So there are constants $r>0$, 
$\epsilon<a$ (depending on $r$) and $\delta<\epsilon/3$ (depending on $\epsilon$) 
such that the following relations hold, where $U_r:=h^{-1}(N_r(0)$ 
\[N_\epsilon(U_r\cap K)\subset h^{-1}(N_{2r}(0))\text{ and }
N_{3\delta}(f_i^{-1}(U_r\cap K)\subset f_i^{-1}(N_\epsilon(U_r\cap K))\]
In the definition of $h_*$ we choose the covering isometry $V$ such that $V^*$ maps 
$L^2(N_r(0))$ into $L^2(U_r)$.  
In the construction of $f_{i*}$ we use an isometry $V_i$ between $L^2(M_i, E_i)$ and $L^2(M)$  
that covers $f_i$ with propagation speed less than $\delta$. 
Let $\chi$ be a normalization function whose Fourier transform is 
supported in $(-\delta,\delta)$. Then $\chi(\mathbf D_i)$ and hence 
$(\alpha(\mathbf D_1), \alpha(\mathbf D_2)\in D^*(M_1,E_1)\times D^*(M_2,E_2)$ has 
propagation speed less than $\delta$. With above assumptions, for  $\xi\in L^2(N_r(0))$ 
the values of $\alpha(\mathbf D_i)(V_i^*(V^*(\xi)_{|K}))$; and hence the value of 
$p(\alpha(\mathbf D_i)(V_i(V(\xi))), \alpha(\mathbf D_i)(V_i(V(\xi))))$ which is 
$h_*([\mathbf D_1, \mathbf D_2])(\xi)$ depends only on the geometric data in  
the neighbourhood $(h\circ f_i)^{-1}(N_{2r}(0))$. This neighbourhoods can
be made arbitrary small by making $r$ as small as necessary. 
\end{pf}

It is clear that by changing the geometric data in a bounded neighbourhood of $N_i$ and $N$
the coarse classes of $f_i$'s and of $h$ do not change and $\mathbf D_i$ change continuously. 
So we assume that the geometric data  around partitioning hyper-surfaces $N_i$'s 
and $N$ take product forms. Therefore in view of above proposition we can assume 
that the whole geometric 
data are of product form coming from partitioning hyper-surfaces. 
More precisely, 
 let $(N,g)$ and $(N_i,g_i)$ for $i=1,2$ be complete Riemannian manifolds with 
 compact subsets $Z\subset N$ and $Z_i\subset N_i$ and smooth maps $f_i\colon (N_i,Z_i)\to (N,Z)$. 
 Moreover let $E_i\to N_i$ be Clifford bundles with corresponding Dirac 
 operators $D_i$. We assume an 
 isometry $\psi:N_1 \verb+\+ Z_1\to N_2 \verb+\+ Z_2$ that is covered by a bundle 
 isometry $\Psi$ 
 such that on smooth sections of $E_1$ which are supported in 
 $N_1 \verb+\+ Z_1$ one has $\Phi^{-1}D_2\Phi=D_1$. 
 Moreover we assume that $f_2\circ \phi=f_1$. Given all these data the relative index 
 $\ind(D_1,D_2)$ is an element in $K_0(C^*(Z\subset N))\simeq\Z$.  
 
 Now let $M=N\times\R$ with product metric $g+(dx)^2$  
 and similarly for $M_i:=N_i\times\R$ and put 
 $W:=Z\times\R$ and $W_i:=\Z_i\times\R$. 
 All geometric structures, including metrics, Clifford bundles and    
 functions $f_i$ can be extended to the product spaces $M_i$. 
 For example $\psi\times\Id:M_1\verb+\+ W_1\to M_2\verb+\+ W_2$ 
 is an isometry which is covered by $\Psi$. 
  Let $\mathbf D_i$ denote the Dirac type operator acting on $C^\infty(M_i,E_i)$. 
 Therefore we have again $\Psi^{-1}\mathbf D_2\Psi=\mathbf D_1$ on 
 $C^\infty(M_1\verb+\+ W_1)$. Given these data, the relative coarse index 
 $\ind(\mathbf D_1,\mathbf D_2)$ is an element in $K_1(C^*(W\subset M))$. 
 The following proposition is a main step toward the complete proof 
 for the main theorem \ref{mainthm}.
 \begin{prop}\label{reduccyl}
 In the above product situation the following relation holds in $K_0(Z\subset N)$
  \[\partial_{mv}(\ind(\mathbf D_1,\mathbf D_2))=\ind(D_1,D_2)\]
 \end{prop}
 \begin{pf}
 The subsets $\R^\pm\times N_i$ and $\R^\pm\times N$ provide respectively excisive partitioning  
 for $M_i=\R\times N_i$ and $M=\R\times N$.  These subsets are flasque, so we are in the situation 
 of corollary \ref{meyerassem}. Therefore, by using this corollary and lemma \ref{asemrel},  it is enough to show 
 the following similar 
 result in $K$-homology  $K_0(Z\subset N)$
 \[\partial_{mv}([\mathbf D_1,\mathbf D_2])=[D_1,D_2]\]
 The decomposition $\R=\R^+\cup \R^-$ provides flasque excisive decomposition for the product 
Manifold $M=\R\times N$. This decomposition provides Mayer-Vietoris exact sequence 
for $\mathcal A$, $\mathcal B$, $\mathcal B/\mathcal A$ and for $\R\times N$ 
and provides following commutative diagram with vertical isomorphisms 
\begin{equation}\label{Calgary}
\begin{CD}
K_0(\mathcal B/\mathcal A) @>>> K_1(\mathcal A) @> p >> K_1(\R\times Z\subset \R\times N)\\
@V{\partial_{mv}}VV @V{\partial_{mv}}VV  @V{\partial_{mv}}VV\\
K_1(B/A) @>>> K_0(A) @> p >> K_0(Z\subset N)
\end{CD}
\end{equation}
Here $p$ is the projection introduced right after \eqref{dirs}.
Using this commutative diagram it is enough to prove the relation 
\begin{equation}\label{sarshar}
\partial_{mv}([\alpha(\mathbf D_1),\alpha(\mathbf D_2)])=
[\beta(D_1), \beta(D_2)]\in K_1(B/A).
\end{equation}
For this purpose we recall the definitions of algebras $\mathcal A$ and $\mathcal B$ just before \eqref{relcnc} 
and the definitions of algebras $A$ and $B$ in the remark after \eqref{dirs}. 
By ignoring the equality on $W$, we have morphisms $j$ from $\mathcal B/\mathcal A$ and $B/A$ 
respectively into $(D^*(M)/C^*(M))\times(D^*(M)/C^*(M)$ and 
$(D^*(N)/C^*(N))\times (D^*(N)/ C^*(N)$ and similarly for other algebras. In particular this is also true for the algebras that 
are involved in the construction of the Mayer-Vietoris morphisms, c.f., \cite[section 3]{siegel2012mayer}. 
It is clear that  $j(\alpha(\mathbf D_1),\alpha(\mathbf D_2))=(\alpha(\mathbf D_1),\alpha(\mathbf D_2))\in 
 (D^*(M)/C^*(M))\times(D^*(M)/C^*(M))$ that represents 
 $([\mathbf D_1], [\mathbf D_2])$ in $K_0(\R\times N)\oplus K_0(\R\times N) $. 
 By \cite[Lemma 4.6]{siegel2012mayer} we have  $\partial_{MV}([\mathbf D_1],[\mathbf D_2])=([D_1], [D_2])\in 
  K_1(D^*(N)/C^*(N))\times (D^*(N)/ C^*(N))$ and this class is represented by $(\beta(D_1),\beta(D_2))$, a projection in 
   $(D^*(N)/C^*(N))\times (D^*(N)/ C^*(N))$. By our discussion on the inclusion of those algebras that are involved in the 
 construction of the Mayer-Vietoris morphism in their product counterparts (by ignoring the equality on some subsets), 
  we conclude that the image of  $[\alpha(\mathbf D_1),\alpha(\mathbf D_2)]\in K_0(\mathcal B/\mathcal A)$, under 
  the Mayer-Vietoris morphism  in the left side of \eqref{Calgary}, is also represented by $((\beta(D_1),\beta(D_2))$, 
  which is the same as  \eqref{sarshar}. 
\end{pf}

Now we have every things to give a complete proof for the main theorem.

\textbf{Proof of the main theorem \ref{mainthm}~:} In proposition \ref{redtopro} we 
showed that $h_*(\ind(\mathbf D_1,\mathbf D_2))$ does not change if we replace data 
with corresponding cylindrical form discussed just before the proposition \ref{reduccyl}. 
Of course this is the case for $\ind(D_1, D_2)$. 
Therefore it is enough to prove the theorem for the product case.  
The following commutative diagram follows from the naturality of the 
Mayer-Vietoris morphisms where $P_*$ is defined in \eqref{conv2}
\begin{equation*}
\begin{CD}
K_1(\R\times Z\subset \R\times N) @>h_*>> K_1(C^*(\R)) \\
@V{\partial_{mv}}VV        @V{\partial_{mv}}VV \\
K_0(Z\subset N) @>P_*>> K_0(C^*(pt.))
\end{CD}
\end{equation*}
Using this diagram and the proposition \ref{reduccyl} we get 
\[\partial_{mv}\circ h_*(\ind(\mathbf D_1 ,\mathbf D_2))=
P_*\circ\partial_{mv}(\ind(\mathbf D_1 ,\mathbf D_2))=P_*(\ind (D_1,D_2))\]
This is precisely the equality \eqref{mainequ} up to isomorphisms \eqref{conv1} and 
\eqref{conv2} and completes the proof of the theorem. 

\begin{rem}
So far we have considered manifolds which are partitioned by one hyper-surface 
(here $N_i$ and $N$). However every things we have proved generalize readily to a 
more general situation, where there are several partitioning hyper-surfaces. 
More precisely for $i=1,2$ let $(M_i, g_i)$, $(M,g)$, $W_i$, $f_i$, $\psi$, $\Psi$, $E_i$, 
and $\mathbf D_i$  be as we have stated in the introduction. 
We assume $n$ be an odd integer and $q$ in bellow be even,  
although these conditions can be relaxed. 
Therefore, as before, we can define the relative index $\ind(\mathbf D_1,\mathbf D_2)$ 
in $K_n(W\subset M)$. 
In what follows we use the subscript $0$ for $M$ and structures on it. For example 
$W_0$ stands for $W$ and $M_0$ stands for $M$ and etc.  
Now let $N_1^1,\dots, N_1^q$ (resp. $N_2^1,\dots, N_2^q$ and $N_0^1,\dots, N_0^q$) be $q$ 
hyper-surfaces in $M_1$ (resp. in $M_2$ and $M_0$) such that $q$ is odd and for $i=0,1,2$: 
\begin{enumerate}
\item $N_i^1,\dots, N_i^q$ intersect each other transversally and coarsely, 
so $N_i=\cap_j N_i^j$ is a sub-manifold of $M_i$ with codimension $q$
\item $N_i^1,\dots, N_i^q$ intersect $W_i$ transversally and coarsely and 
$Z_i=W_i\cap(\cap_q N_i^q)$ is compact, 
\item for $i=1,2$ and $j=1,\dots, q$ we have $N_i^j=f_i^{-1}(N_0^j)$ and $N_0^j$ are 
regular for $f_i$. 
\end{enumerate}
For $i=0, 1,2$ the sub manifolds $N_i$ are equipped with the induced Clifford bundles 
and have their Dirac operators $D_i$ which are conjugate outside $Z_i$. Therefore 
we can define the relative index $\ind(D_1,D_2)$ as an element in $K_p(C^*(Z\subset N))$, where 
$p=n-q$. 
With above assumptions the signed distances from $N_0^1, \dots, N_0^q$ define the 
map $h:M\to \R^q$ which is coarse outside $W=W_0$. Therefore it induces a map 
$h_*:K_*(C^*(W\subset M))\to K_n(C^*(\R^q))\simeq \Z$. 
 Under these conditions we have 
\begin{equation}\label{genmain}
h_*(\ind(\mathbf D_1, \mathbf D_2))=\ind(D_1,D_2)
\end{equation}
This is a generalization of \eqref{mainequ} and its proof is completely similar.  
\end{rem}

\section{Application to positive scalar curvature problem}\label{sec.5}

In this section we use the theorem \ref{mainthm} to prove that $\tilde N\times \R$ can not have a 
uniformly positive scalar curvature if $\tilde N$ is an enlargeable manifold. 
This theorem was stated and proved by means of geometric tools 
(and the original version of relative index theorem) in \cite{GrLa3}. 
It was also reproved in \cite{Zadeh-Index} and \cite{Zadeh-Note} by using a 
version of the partitioned index theorem 
with general $C^*$-algebra coefficients. To provide yet another proof for this fact we state and 
prove an expected vanishing theorem concerning the relative index. 
As $\mathbf D_i$ is a Dirac type operator, due to the Weitzenbock formula
\[\mathbf D_i^2=\nabla\nabla^*+R_i~,\]
where $R_i$ is the Clifford curvature; a self-adjoint bundle linear map on $E_i$. If 
$E_i$ is the spin bundle associated to a spin structure, or a 
spin bundle twisted by a flat bundle, 
then $R_i=\kappa_i/4$, where $\kappa_i$ is the scalar curvature of the underlying 
Riemannian manifold $(M_i,g_i)$. 
\begin{thm}\label{vanish}
If for $i=1,2$ the Clifford curvature $R_i$ are uniformly positive on $M_i$ then 
\[\ind(\mathbf D_1,\mathbf D_2)=0\in K_1(C^*(W\subset M))~.\] 
\end{thm}
\begin{pf}
The assumption on the Clifford curvature along with the Weitzenbock formula 
implies that there is a gap around $0$ in the spectrum of $\mathbf D_i$. So, the 
normalizing function $\chi$ in the definition of the relative index 
(just before relation \eqref{relcnc}) can be assumed to satisfy 
$\chi^2=1$. Therefore $(\alpha(\mathbf D_1), \alpha(\mathbf D_2))$ itself is a projection in 
$\mathcal B$ and the class $[(\alpha(\mathbf D_1), \alpha(\mathbf D_2))]$ in 
$K_0(\mathcal B/\mathcal A)$ is the image of a class in $K_0(\mathcal A)$. 
This implies the vanishing of $\delta([(\alpha(\mathbf D_1), \alpha(\mathbf D_2))])$ in 
\eqref{relcnc} and then the vanishing of the relative index. 
\end{pf}

Following Gromov-Lawson \cite{GrLa3} ,
let $\tilde N$ be a closed oriented manifold of dimension $n$ with  
a fixed riemannian metric $\tilde g$. The manifold $\tilde N$ is enlargeable if for each 
real number $\epsilon>0$ there is a Riemannian spin cover $(N, g)$, 
with lifted metric, and a smooth map $f:N\to S^n$ such that: 
the function $f$ is constant outside a compact subset $Z$ of $N$; 
the degree of $f$ is non-zero;  and
the map $f:(N, g)\to(S^n,g_0)$ is $\epsilon$-contracting, where
$g_0$  is the standard metric on $S^n$. Being $\epsilon$-contracting means that 
$\|T_xf\|\leq\epsilon$ for each $x\in N$, where $T_xf\colon T_x N\to T_{f(x)}S^n$. 

With above notation, it turns out that  there is 
a Hermitian bundle $E\to  N$ which is isomorphic 
to the trivial bundle $F:=C^k\times  N$  
outside the compact subset $ Z\subset N$ such that $\omega:=\Ch(E-\C^k)$ is 
a non-zero multiple of the volume element of $N$ which vanishes outside $ Z$. 
Moreover the bundle $E$ has a connection 
whose curvature $R$ is bounded from above by $\epsilon$, i.e. $\|R\|\leq\epsilon$.

\begin{thm}\label{theo2}
For an enlargeable closed manifold $\tilde N$, the product space $M=\tilde N\times \R$ 
does not admit a Riemannian metric with uniformly positive scalar curvature in the coarse equivalent 
class of $dt^2+g_{\tilde N}$. 
\end{thm}
\begin{pf}
Let $D_1$ and $D_2$ stand 
for Dirac operators of $N$ twisted; respectively; by $E$ and $F$. In this situation 
the relative index $\ind_r(D_2,D_1)\in K_0(Z)\simeq\Z$ is given by the following 
relation \cite[formula 4.5]{Roe-relative} 

\begin{equation}\label{nonvan}
 \ind_r(D^E,D^F)=\int_{\tilde N}A(TM)\wedge \Ch(E-F)=\int_{\tilde N}\omega\neq0
\end{equation}

Now consider the manifold $M=N\times \R$ and put $W:=Z\times R$. 
Because $Z$ is compact, with respect to any lifted (from $\tilde N\times \R$) Riemannian metric on $M$ the intersection 
of $N$ and $Z\times\R$ is transverse and coarse.
Bundles $E$ and 
$F$ will be considered as bundles over $M$. Then the data $(M,W,\mathbf D_1)$ 
and $(M,W,\mathbf D_2)$ 
satisfy the relative index theorem conditions with $f_1=f_2=\Id$ and 
$h:N\times\R\to\R$ being the projection on the second factor. The theorem \ref{mainthm} reads 
\[h_*(\ind(\mathbf D_1,\mathbf D_2))=\ind(D_1,D_2)\neq 0\]  
However if the scalar curvature of a metric on $\R\times \tilde N$ is 
everywhere uniformly positive, then 
it is so for lifted metric on  $M=N\times \R$. 
The  curvature of $E\to M$ is bonded by $\epsilon$. 
Therefore if $\epsilon$ is sufficiently small, the Weitzenbock formula and lemma 
\eqref{vanish} together provide the vanishing result  $\ind(\mathbf D_1,\mathbf D_2)=0$. 
This is in contradiction with \eqref{nonvan} and implies the non-existence of such a metric. 
\end{pf} 

\begin{rem}
Using the equality \eqref{genmain}, we can prove the following generalization of 
the above theorem by a similar proof: Let $\tilde N$ be an enlargeable manifold and consider 
$M=\tilde N\times \R^q$. Let $\tilde g$ be a riemaninan metric on $\tilde N$ 
and $g$ be a riemannian metric on $M$ which is in the same 
coarse class that a product metric $\tilde g+dx_1^2+\dots+ dx_q^2$ 
(these later metrics all fall into the same coarse class whatever $\tilde g$ might be). Then 
the $q$ hyper-surfaces $N_i=N\times \R$ intersects coarsely and we can, as in above, 
apply \eqref{genmain} to conclude that the scalar curvature of $g$ cannot be uniformly positive. 
Despite this fact, it is proved in \cite{GrLa3} that $M$ does admit a metric 
with uniformly positive scalar curvature, provided that $q\geq 3$. 
\emph{Consequently this metric is not in the same coarse class as the product 
metric $\tilde g+dx_1^2+\dots+ dx_q^2$}. 
\end{rem} 


\end{document}